\documentclass[11pt]{article}
\usepackage{amsfonts}
\usepackage{mathrsfs}
\usepackage{amsmath,amssymb, float}
\usepackage{mathtools}

\DeclareMathOperator{\EX}{\mathbb{E}}
\DeclareMathOperator{\E}{\mathbb{E}}
\usepackage[latin1]{inputenc}
\usepackage{amsmath,amssymb}
\usepackage{latexsym}
\usepackage{epstopdf}
\usepackage{geometry}   
\usepackage[active]{srcltx}
\usepackage{graphicx}
\usepackage{epstopdf}
\usepackage[
bookmarks=true,         
bookmarksnumbered=true, 
colorlinks=true, pdfstartview=FitV, linkcolor=blue, citecolor=blue,
urlcolor=blue]{hyperref}

\newtheorem {theorem}{Theorem}[section]

\newtheorem{assumption}{Assumption}
\newtheorem {corollary}{Corollary}[section]

\newtheorem{lemma}{Lemma}[section]

\newtheorem{remark}{Remark}[section]

\newenvironment{proof}[1][Proof]{\textbf{#1.} }{\
\rule{0.5em}{0.5em}}

\def\R{{\mathbb R}}
\def\E{{{\mathbb E}\,}}
\def\P{{\mathbb P}}

\def\Z{{\mathbb Z}}

\def\N{{\mathbb N}}
\def\Var{{\mathop {{\rm Var\, }}}}

\date{}
\begin{document}

\title{On the integrated mean squared error of wavelet density estimation for linear processes}

\maketitle
\begin{center}
\bigskip Aleksandr Beknazaryan$^{a}$, Hailin Sang$^{b}$ and Peter Adamic$^{c}$

\bigskip$^{a}$ Institute of Environmental and Agricultural Biology (X-BIO), University of Tyumen, Tyumen, Russia,
a.beknazaryan@utmn.ru

\bigskip$^{b}$ Department of Mathematics, University of Mississippi, University, MS 38677,  USA
sang@olemiss.edu

\bigskip$^{c}$   Department of Mathematics and Computer Science, Laurentian University, Sudbury, Ontario, Canada, padamic@cs.laurentian.ca
 \bigskip
\end{center}

\begin{abstract}
\noindent  Let $\{X_n:  n\in \N\}$ be a linear process with density function $f(x)\in L^2(\R)$. We study wavelet density estimation of $f(x)$. Under some regular conditions on the characteristic function of innovations, we achieve, based on the number of nonzero coefficients in the linear process, the minimax optimal convergence rate of the integrated mean squared error of density estimation. Considered wavelets have compact support and are twice continuously differentiable. The number of vanishing moments of mother wavelet is proportional to the number of nonzero coefficients in the linear process and to the rate of decay of characteristic function of innovations. Theoretical results are illustrated by simulation studies with innovations following Gaussian, Cauchy and chi-squared distributions.
\vskip.2cm

\vskip.2cm \noindent {\bf Keywords}:
\noindent  linear process, wavelet method, density estimation, projection operator


\vskip.2cm 

\end{abstract}

\section{Introduction}\label{introduction}
In this paper we consider the linear process 
\begin{equation}\label{lp}
X_n=\sum\limits^{\infty}_{i=0} a_i\varepsilon_{n-i},
\end{equation} 
with probability density function $f(x)$, where the innovations $\varepsilon_i$ are independent and identically distributed (i.i.d.) real-valued random
variables in some probability space $(\Omega, \mathcal{F}, \P)$ and $\{a_i\}$ are real coefficients such that $a_0\ne 0$ and $\sum\limits_{i=0}^\infty a_i \varepsilon_{n-i}$ is well defined. With observations $\{X_k\}_{k=1}^n$, we introduce a wavelet based estimator of density function $f(x)$ and study the order of its integrated mean squared error (IMSE). 

Most of the estimators for such processes presented in the literature are the kernel estimators. In particular, under the condition that the innovations have bounded second moment,  central limit theorems for kernel density estimators of linear processes were derived in Wu and Mielniczuk (2002). Hall and Hart (1990) proved that the IMSE of the kernel estimator for two-sided linear processes with given $n$ observations can be decomposed into the sum of $2$ terms: the IMSE of the same kernel estimator based on i.i.d. sample of size $n$ following the same distribution and a term proportional to the variance of the sample mean. Hence, since the variance of the sample mean of short range dependent linear processes has the order $O(1/n)$ (see, e.g., Priestley (1981)), the IMSE of the kernel estimator for two-sided short range dependent linear processes has the same order as the IMSE of the kernel estimator for i.i.d. data. Analogous results for a larger class of processes are given by Mielnizcuk (1997), where under certain conditions on the Fourier transform of the kernel estimator, a similar decomposition for the IMSE is obtained. For strictly stationary sequences, the rates of the IMSE of kernel density estimators are derived in Meloche (1990). Moving average processes
$$z_t=\varepsilon_t-\theta_1\varepsilon_{t-1}-...-\theta_q\varepsilon_{t-q}$$ 
of order $q$ $(\ge 1)$ has been studied in Saavedra and Cao (2000). Assuming that the density function of the innovation $\varepsilon_1$ has continuous and bounded derivatives of order up to $4$, the parametric order $1/n$ for the IMSE of kernel density estimator was derived. 

For independent sequences, the IMSE rates of wavelet density estimators  have been studied in the works Donoho and  Johnstone (1994, 1995), Donoho et al. (1995, 1996), Zhang and Zheng (1999), Lu (2013), and Gin\'{e} and Madych (2014). However, compared to the kernel density estimators, for linear processes the IMSE rates of wavelet density estimators have been barely covered in the literature. Below we present some of the works where related topics are studied. First, a general overview of the use of wavelets in the theory of non-parametric function estimation is given in Johnstone (1999). For time series data, wavelet-based estimators of densities and marginal densities were obtained in Leblanc (1996) and Gannaz and Wintenberger (2010). Note that the random variables considered in those two works are assumed to have compactly supported distributions. Chesneau (2014) studied the wavelet hard thresholding density estimators for observations with $\alpha$-mixing dependence. Also, under some regular conditions on the strong mixing coefficient and on the densities, Kou and Guo (2018) gave an upper bound 
for the IMSE of linear wavelet density estimators of strong mixing sequences. In particular, if $a_k=2^{-k}$ and $\varepsilon_k\sim N(0,\sigma^2),$ $k\geq 0,$ then the linear process  \eqref{lp} is strong mixing. Wavelet based methods for various mixing, quadrant and $m$-dependence conditions have also been developed in the works Badaoui and Rhomari (2015), Chesneau (2012) and Chesneau, Dewan and Doosti (2012),  Chesneau, Doosti and Stone (2019) and Li and Zhang (2022). In particular, in the recent related works Chesneau, Doosti and Stone (2019) and Li and Zhang (2022), wavelet-based estimators of bounded, compactly supported densities of $m$-dependent processes and of spectral densities of non-Gaussian linear processes belonging to Besov spaces are proposed. The rates of convergence derived for those estimators coincide with the minimax optimal rates with additional logarithmic factors.

In the present work we study the IMSE of wavelet density estimation for linear processes of the form \eqref{lp}. The estimator is obtained by substituting the coefficients in wavelet series expansion of density function $f\in L^2(\mathbb{R})$ by the average values of corresponding translations and dyadic
dilations of wavelets at the observations  $\{X_k\}_{k=1}^n$. Optimal choice of wavelets is determined by the rate of decay of the characteristic function of innovations $\varepsilon_i$ and by the number of nonzero
coefficients in the linear process. Note that those two features of the linear process also determine the order of smoothness of the density function $f$. In particular, if the induced smoothness of $f$ is of order $m$, then the optimal choice of wavelets allows to achieve the minimax optimal nonparametric estimation rate $n^{-\frac{2m}{2m+1}}$. In case $m=\infty$, the IMSE rates of estimators can get arbitrarily close to $n^{-1}$, and in this case the wavelets with higher number of vanishing moments induce better rates of estimation.
Note that there is a similar phenomena for kernel density estimations: kernel functions with higher order allow to derive faster estimation rates (see, e.g., Marron (1994)).

The paper is structured as follows: in Section \ref{mainresult} we present the general setting and the main results. We perform simulation studies in Section \ref{simulation} to confirm the results of Section  \ref{mainresult}.  Section \ref{proof} gives the proofs  of several auxiliary results which lead to the proof of the main theorem. The proof is based on the decomposition of the IMSE of linear wavelet estimator into three parts, first of which corresponds to the scaling function and the other two are the partial sum and the tail of the series with coefficients corresponding to mother wavelet. Convergence rates of these parts are derived in Section \ref{proof} by rewriting them as sums of Fourier transforms at integers and applying the formula of Poisson. Section \ref{proof} also presents the properties of considered wavelets, the justification of application of Poisson summation formula and the proofs of several properties of characteristic functions of $X_1$ and $\varepsilon_1$.

We apply the following notations throughout the paper: for a function $g(x)\in L^1(\R)$, we use $\widehat{g}(u)=\int_{\R} e^{-\iota x u}\, g(x)\, dx$ as the definition of  its Fourier transform, where $\iota$ is the imaginary unit. The characteristic function of linear process $X_n=\sum\limits^{\infty}_{i=0} a_i\varepsilon_{n-i}$ is denoted by $\phi(\lambda),$ and $\phi_{\varepsilon}(\lambda)$ denotes the characteristic function of innovations $\varepsilon_i$. That is, $\phi(\lambda)=\E[e^{\iota \lambda  X_1}]$ and $\phi_{\varepsilon}(\lambda)=\E[e^{\iota \lambda  \varepsilon_1}]$.   For each $i, j\in\N$, we define $H(X_i):=H(X_i)(\lambda):=e^{\iota \lambda X_i}-\phi(\lambda)$,   and 
$H_{ij}(u,v):=\EX e^{\iota uX_i+\iota vX_j}-\phi(u)\phi(v)=\EX \big[H(X_i)(u)H(X_j)(v)\big ]$, 
$\lambda, u, v\in\R$.
For two sequences $\{a_n\}$ and $\{b_n\}$, the notation $a_n=O(b_n)$ indicates the existence of a constant $C$ such that $a_n\leq Cb_n$ for all $n\in\mathbb{N}$. Finally, $\lceil \cdot\rceil$ and $\lfloor \cdot\rfloor$ denote, respectively, the usual ceiling and floor functions on $\mathbb{R}$, and the positive constant $c$ that appears in the proofs may vary from line to line, but is independent of $i, j, k, n$ and other indices involved.

\section{Main results}\label{mainresult}
Before proceeding to the imposed assumptions and the main result, let us present several particular settings for which the linear process \eqref{lp} is well-defined. If $\{\varepsilon_{i}:\; i\in\Z\}$ is a sequence of i.i.d. random variables in $L^{p}(\R)$ for some $p>0$, $\E\varepsilon_{i}=0$ when $p\geq 1$, and
$\{a_{i}\}^{\infty}_{i=0}$ is a sequence of real coefficients such that $\sum\limits_{i=0}^{\infty}|a_{i}
|^{2\wedge p}<\infty$, by Kolmogorov's three-series theorem,  the linear process $X_n$ given in (\ref{lp}) 
exists and is well-defined. 
As for linear processes with symmetric $\alpha$-stable innovations $(0<\alpha< 2)$, i.e., the law of innovations having characteristic function $\E[e^{\iota\lambda\varepsilon_1}]=\exp(-c_{\alpha}|\lambda|^\alpha)$ for some positive constant $c_{\alpha}$ only depending on $\alpha$, by Kolmogorov's three-series theorem, the series given in \eqref{lp} converges almost surely if and only if $\sum\limits_{i=0}^{\infty}|a_{i}
|^{\alpha}<\infty$ (see, e.g., Samorodnitsky and Taqqu (1994)).

We now develop the wavelet estimator $\hat{f}_n$ of the density function $f$ of linear process \eqref{lp}. As the density function $f$ is assumed to be in $L^2(\mathbb{R}),$ then $f$ accepts the representation 
\begin{equation}\label{1}
f(x)=\sum\limits_{k\in\mathbb{Z}}\alpha_{0k}\varphi_{0k}(x)+\sum\limits_{j\geq 0}\sum\limits_{k\in\mathbb{Z}}\beta_{jk}\psi_{jk}(x),
\end{equation}
where $\psi_{jk}(x)=2^{j/2}\psi(2^jx-k), \varphi_{jk}(x)=2^{j/2}\varphi(2^jx-k), j, k\in \mathbb{Z},$ are the translations and dyadic dilations of mother wavelet $\psi$ and the scaling function $\varphi$. The orthonormality of functions $\varphi_{0k}(x)$ and $\psi_{jk}(x)$, $j\ge 0, k\in \Z$, in \eqref{1} implies that
\begin{equation*}
\alpha_{0k}=\int_\mathbb{R} f(x)\varphi_{0k}(x)dx=\EX \varphi_{0k}(X_1),
\end{equation*}
\begin{equation*}
\beta_{jk}=\int_\mathbb{R} f(x)\psi_{jk}(x)dx=\EX \psi_{jk}(X_1).
\end{equation*}
Therefore, given a sample $X_1,...,X_n,$ we can approximate the coefficients $\alpha_{0k}$ and $\beta_{jk}$ by the quantities $\hat{\alpha}_{0k}$ and $\hat{\beta}_{jk}$ defined as
\begin{equation*}
\hat{\alpha}_{0k}=\frac{1}{n}\sum\limits_{i=1}^n\varphi_{0k}(X_i)
\end{equation*}
and
\begin{equation*}
\hat{\beta}_{jk}=\frac{1}{n}\sum\limits_{i=1}^n\psi_{jk}(X_i), 
\end{equation*}
and obtain an estimator  
\begin{equation}\label{estimator}
\hat{f}_n(x)=\sum\limits_{k\in\mathbb{Z}}\hat{\alpha}_{0k}\varphi_{0k}(x)+\sum\limits_{j= 0}^{j_n}\sum\limits_{k\in\mathbb{Z}}\hat{\beta}_{jk}\psi_{jk}(x)
\end{equation}
of $f(x)$. The choice of considered wavelets is specified in the Assumption \ref{assu}, while the selection of optimal truncation location $j_n$ is given in the Theorem \ref{main} below. Assumptions on the linear process (\ref{lp}) needed to derive the IMSE rate of the estimator (\ref{estimator}) are presented in Assumption \ref{assu0}.

\begin{assumption}\label{assu0}
For the linear process \eqref{lp}, suppose that the density function  $f(x)$ is in $L^2(\R)$ and the coefficients satisfy $\sum\limits^{\infty}_{i=0}|a_i|^{\gamma}<\infty$, for some $\gamma \in (0,1]$. Suppose that there are at least $M\geq 1$ non-zero coefficients, and, without loss of generality, let $a_0\ne 0$. Assume that for some constant $\beta\geq\gamma$ the characteristic function $\phi_{\varepsilon}$ of innovations satisfies     
	\begin{align}\label{integral}
 u^{\beta}\phi_{\varepsilon}(u)\in L^1(\mathbb{R})\cap L^\infty(\mathbb{R}).
	\end{align}

\end{assumption}

\begin{assumption}\label{assu}
We consider wavelets for which both $\psi$ and $\varphi$ have compact support and are twice continuously differentiable. Also, we assume that $\psi$ has $\lceil{M\beta}\rceil$ vanishing moments, where $M$ and $\beta$ are the parameters from Assumption \ref{assu0}.
\end{assumption}
Note that the Daubechies wavelets of order $2\lceil{M\beta}\rceil\lor 8$ and higher satisfy the Assumption \ref{assu} (see Daubechies (1992), Chap. 7).

The following theorem is the main result of the paper.
\begin{theorem} \label{main} 
 Assume that the derivative of the characteristic function $\phi_{\varepsilon}$ of innovations has bounded modulus on $\mathbb{R}$.  
In case  $\gamma\in (1/2,1]$, assume that the innovations also satisfy  
\begin{align}\label{cha2}
\E|e^{\iota \lambda \varepsilon_1}-\phi_{\varepsilon}(\lambda)|^{2}\leq c \left(|\lambda|^{2\gamma}\wedge 1\right)
\end{align}
for all $\lambda\in \R$ and for some constant $c>0$ which only depends on $\gamma$. 
Then, under the Assumptions  \ref{assu0} and \ref{assu},  for $j_n=\lceil\frac{\log_2n}{2M\beta+1}\rceil$,
	\begin{align*} 
	\EX \int_\mathbb{R} \bigg[\hat{f}_n(x)-f(x)\bigg]^2dx=O\bigg(n^{-\frac{2M\beta}{2M\beta+1}}\bigg).
	\end{align*}
\end{theorem}

Note that the derivative of $\phi_{\varepsilon}$ has  bounded modulus on $\mathbb{R}$ if, for example, $\EX |\varepsilon_i|<\infty$ (see Lukacs (1996), page 22). Also, the condition \eqref{cha2} is satisfied if, for example, $\E |\varepsilon|^{2\gamma}<\infty$ (see Lemma 7.3 in Sang, Sang and Xu (2018)). We thus have the following corollary of Theorem \ref{main}:
\begin{corollary}Assume that the innovations $\varepsilon_i$ of the linear process \eqref{lp} have finite second moment. Then, under the Assumptions  \ref{assu0} and \ref{assu},  for $j_n=\lceil\frac{\log_2n}{2M\beta+1}\rceil$,
	\begin{align*} 
	\EX \int_\mathbb{R} \bigg[\hat{f}_n(x)-f(x)\bigg]^2dx= O\bigg(n^{-\frac{2M\beta}{2M\beta+1}}\bigg).
	\end{align*}
\end{corollary}


\begin{remark}\normalfont{If the linear processes have finitely many non-zero coefficients, say, $M$, then the condition \eqref{integral} implies
\begin{equation}\label{charineq}
\int_\mathbb{R}|u^{M\beta}\phi(u)|du<\infty,
\end{equation} 
where $\phi$ is the characteristic function of the linear process $X_i$. \eqref{charineq}, in turn, implies (see Folland (1999), Theorem 8.22) that $f\in C^{M\beta}$. Remarkably, the rate $n^{-\frac{2M\beta}{2M\beta+1}}$, which is obtained by applying Daubechies wavelets of order $2\lceil{M\beta}\rceil\lor 8,$ is the best possible mean square convergence rate to a $C^{M\beta}$-smooth density even with given $n$ \textit{independent} observations (see, e.g., Wahba (1975)). If the linear processes have infinite number of non-zero coefficients, theoretically we can have IMSE rate arbitrarily close to $n^{-1}$. However, we shall apply wavelets with $m$ vanishing moments to reach the rate $n^{-\frac{2m}{2m+1}},$ and the larger $m$ is chosen the closer we get to the optimal rate $n^{-1}$. Note that in kernel density estimation we have similar requirement on the kernel function: to attain the rate $n^{-\frac{2m}{2m+1}}$ of estimation of density function $f\in C^m$, where $m$ is an even positive integer, we shall apply $m$-th order kernel function $K(\cdot)$ with $\int_\mathbb{R}x^pK(x)dx=0$, $1\le p<m$, (see, e.g., Marron (1994)). }
\end{remark}
\begin{remark}
\normalfont{Note that in Theorem \ref{main} the estimator $\hat{f}_n$ depends on the value of $j_n$ which, in turn, is determined by the values of $M$ and $\beta$. Hence, the identification of optimal estimator assumes the knowledge of both $M$ and $\beta$. 
}
\end{remark}
\begin{remark}
\normalfont{Note that if $\gamma\in (0,1/2]$ then the condition \eqref{cha2} is automatically satisfied (see Lemma \ref{case}). If the innovations have non-degenerate distribution, the range $\gamma\in (0,1]$ in inequality  \eqref{cha2} is optimal and cannot be extended (see Lemma 7.4 in Sang, Sang and Xu (2018)). Also, it is easy to see that $$\Var(e^{\iota  \lambda \varepsilon_1})=\E|e^{\iota \lambda \varepsilon_1}-\phi_{\varepsilon}(\lambda)|^{2}=1-|\phi_{\varepsilon}(\lambda)|^{2}.$$
Therefore,  condition \eqref{cha2} together with $\sum\limits^{\infty}_{i=1}|a_i|^\gamma<\infty$ imply \[
\sum\limits^{\infty}_{i=1} \sqrt{\Var(e^{\iota  \lambda a_i \varepsilon_1})}=\sum\limits^{\infty}_{i=1} \sqrt{\E|e^{\iota a_i\lambda \varepsilon_1}-\phi_{\varepsilon}(a_i\lambda)|^{2}}\le c|\lambda|^\gamma\sum\limits^{\infty}_{i=1}|a_i|^\gamma<\infty.
\]
Hence, according to the definition of short or long memory dependence for linear processes as in Sang, Sang and Xu (2018), Theorem \ref{main} provides an IMSE result for short memory linear processes.}
\end{remark}

To see some examples, it is easy to verify that for any $\beta>0$ the condition \eqref{integral} is satisfied for linear processes with innovations following either Gaussian distribution ($\alpha=2$) or symmetric stable distribution with index $1\le\alpha<2$. In particular, for the Gaussian case we have $\E|e^{\iota \lambda \varepsilon_1}-\phi_{\varepsilon}(\lambda)|^{2}\leq c \left(|\lambda|^{2}\wedge 1\right)$, and, therefore, we can take $\gamma=1$ in condition \eqref{cha2}. In the $\alpha$-stable case, $\E|e^{\iota \lambda \varepsilon_1}-\phi_{\varepsilon}(\lambda)|^{2}\leq c_{\alpha} \left(|\lambda|^{\alpha}\wedge 1\right)$, and the condition \eqref{cha2} is satisfied for  $\gamma=\alpha/2$. Also, the derivative of the characteristic function $\phi_{\varepsilon}(u)=e^{-c|u|^\alpha}$ ($1<\alpha\le 2$) exists and is bounded on $\mathbb{R}$. In the case of $1$-stable distribution, the  derivative of the characteristic function $\phi_{\varepsilon}(u)=e^{-c|u|}$ exists and is bounded everywhere except at $0$. However, it still has bounded left and right derivatives at $0$ which is sufficient for the proof of Theorem \ref{main} for this case. If the linear processes have only finite number of nonzero coefficients, then the condition \eqref{integral} of the theorem is in particular satisfied for linear processes with innovations following chi-squared distribution $\chi^2_k$ with $k>2(\beta+1)$ degrees of freedom and for Gamma distribution $\Gamma(k, \theta)$ with shape $k>\beta+1$. In both cases the condition \eqref{cha2} is satisfied for all $\gamma\in (0,1/2]$ (see Lemma \ref{case}).
 

\section{Simulation study}\label{simulation}

In this section we conduct a simulation study to evaluate the  wavelet density estimator of linear processes with various coefficient sequences and with innovations following stable, Gaussian or chi-squared distributions.  For each case, we shall compare the true density function with our wavelet density estimator. 
We first perform simulation study for the linear processes with infinitely many non-zero coefficients and with innovations following stable or Gaussian distributions.  It has two steps: 

{\bf Step 1}: Calculate or approximate the true density function $f(x)$ of the linear process. Here we use the same method as in Sang, Sang and Xu (2018) to find the density function $f(x)$ of the linear process $X_n$ with various coefficient sequences $\{a_i\}$ and with various innovations $\varepsilon_i$. For convenience, we sketch the method below. 

The linear process in (\ref{lp}) with density function $f(x)$ has characteristic function $\phi(t)=\mathbb{E}[e^{\iota tX_1}]$ of the form 
\begin{align*}
\phi(t)=\mathbb{E}[e^{\iota tX_n}]=\prod_{i=0}^\infty\mathbb{E}[e^{\iota t a_i \varepsilon_{n-i}}]=\prod_{i=0}^\infty\phi_\varepsilon(a_it). 
\end{align*}
For any $d<\frac{1}{2}, d\notin\mathbb{Z}$, let  $a_i=\frac{\Gamma(i+d)}{\Gamma(d)\Gamma(i+1)}$, $i\ge 0$. By Stirling's formula, $\Gamma(x)\sim \sqrt{2\pi}e^{-x+1}(x-1)^{x-1/2}$ as $x\rightarrow \infty$, and $a_i\sim i^{d-1}/\Gamma(d)$ as $i\rightarrow \infty$. 
%
%
We perform simulation study for linear processes with Gaussian or Cauchy innovations and with coefficients determined by certain selected values of $d$. 

{\it Case 1}. Suppose the innovations $\{\varepsilon_i\}$ are i.i.d. standard normal random variables. Then, the characteristic function of the innovation is $\phi_\varepsilon(t)=e^{-\frac{t^2}{2}}$ and, therefore,
$\phi(t)=e^{-\frac{t^2}{2}\sum_{i=0}^\infty a_i^2}$.  In this case the linear process $\{X_n\}$ is normally distributed with mean $0,$ variance $A^2:=\sum_{i=0}^\infty a_i^2$ and has density function $f(x)= \frac{1}{\sqrt{2\pi}A}e^{-\frac{x^2}{2A^2}}$. By the Gauss's theorem for hypergeometric series [Gauss (1866)],  
$A^2=\sum_{i=0}^{\infty} a^2_i=\frac{\Gamma(1-2d)}{\Gamma^2(1-d)}$ for any $d<\frac{1}{2}, d\notin\mathbb{Z}$ (see also Bailey (1935) or its direct calculation in Sang and Sang (2017)). In particular, if $d=-0.5$ then $A^2=4/\pi,$ and $A^2=3.39531$ if $d=-1.5$. 
Hence, for $d=-0.5$ \begin{align*}
f(x)=\frac{\sqrt{2}}{4}e^{-\frac{\pi x^2}{8}}
\end{align*}
and for $d=-1.5$ we have
\begin{align*}
f(x)=0.216506 e^{-0.147262 x^2}.
\end{align*} 

{\it Case 2}. If the innovations $\{\varepsilon_i\}$ have i.i.d. symmetric $\alpha$-stable distribution with $0<\alpha<2$, then $\phi_\varepsilon(t)=e^{-|t|^\alpha}$ and 
$\phi(t)=e^{-|t|^\alpha\sum_{i=0}^\infty |a_i|^\alpha}$.
In particular, if $\alpha=1$, we have the standard Cauchy distribution with $\phi_\varepsilon(t)=e^{-|t|}$ and 
$\phi(t)=e^{-|t|\sum_{i=0}^\infty |a_i|}$. In this case, with the help of the software {\it Mathematica}, we can approximate $\sum_{i=0}^\infty |a_i|$ by  $\sum_{i=0}^{100,000} |a_i|=1.99822,$ if $d=-0.5,$ and by $\sum_{i=0}^{100,000} |a_i|=3,$ if $d=-1.5$. Then, the density of the linear process corresponding to the choice $d=-0.5$ is given by 
\begin{align*}
f(x)=\frac{1}{1.99822\pi(1 + (x/1.99822)^2)},
\end{align*}
and for the case  $d=-1.5$ we have
\begin{align*}
f(x)=\frac{3}{\pi(9+x^2)}.
\end{align*}
As in this case $\E|e^{\iota \lambda \varepsilon_1}-\phi_{\varepsilon}(\lambda)|^{2}\leq c_{\alpha} \left(|\lambda|^{\alpha}\wedge 1\right)$, then we can choose $\gamma=\alpha/2$ to fulfill the condition (\ref{cha2}). As $a_i=\frac{\Gamma(i+d)}{\Gamma(d)\Gamma(i+1)}\sim i^{d-1}/\Gamma(d)$, then, to guarantee the convergence of the series $\sum_{i=0}^\infty |a_i|^\gamma$, we should have $d<1-\frac{1}{\gamma}=1-\frac{2}{\alpha}$. In particular, for innovations following Cauchy distribution we have $\alpha=1$, and, therefore, the convergence of the series is guaranteed for memory parameters $d<-1$. In the following step we will see that in the case $d<-1$ the performance of the estimator is indeed better than in the case when $d>-1$.

{\bf Step 2}: In the previous step we considered linear processes with Gaussian and Cauchy innovations and with memory parameters $d=-0.5$ and $d=-1.5$. By applying a modification of the MATLAB code from Fa\"{y} et al. (2009), for each of these cases we produce a linear process $\{X_i\}_{i=1}^n$ with $n=2^{16}$ observations.  For each generated linear process, we apply the $1-$D wavelet estimator from the Wavelet Analyzer App of MATLAB to estimate the true density function.  As required in Assumption \ref{assu}, we use  the Daubechies wavelets of order $8$ and we do not apply any soft or hard thresholding method in the procedure to obtain the estimator (\ref{estimator}). Figures \ref{fig:1} and \ref{fig:2} below show the performance of wavelet density estimation of linear processes with Gaussian and Cauchy innovations for the cases $d=-0.5$ and $d=-1.5$. They confirm the result of Theorem \ref{main} for the cases when innovations follow Gaussian distribution and $d=-0.5, -1.5$ and when the innovations follow Cauchy distribution and $d= -1.5$. When the innovations follow Cauchy distribution and $d= -0.5$, the performance of the estimator is worse than the other three cases. However, since $d=-0.5>-1$, then, as indicated in the last part of previous step, the result in this case is not guaranteed by Theorem \ref{main}. 
\begin{figure}[H]
\centering
\begin{tabular}{cc}
\includegraphics[width=0.45\linewidth]{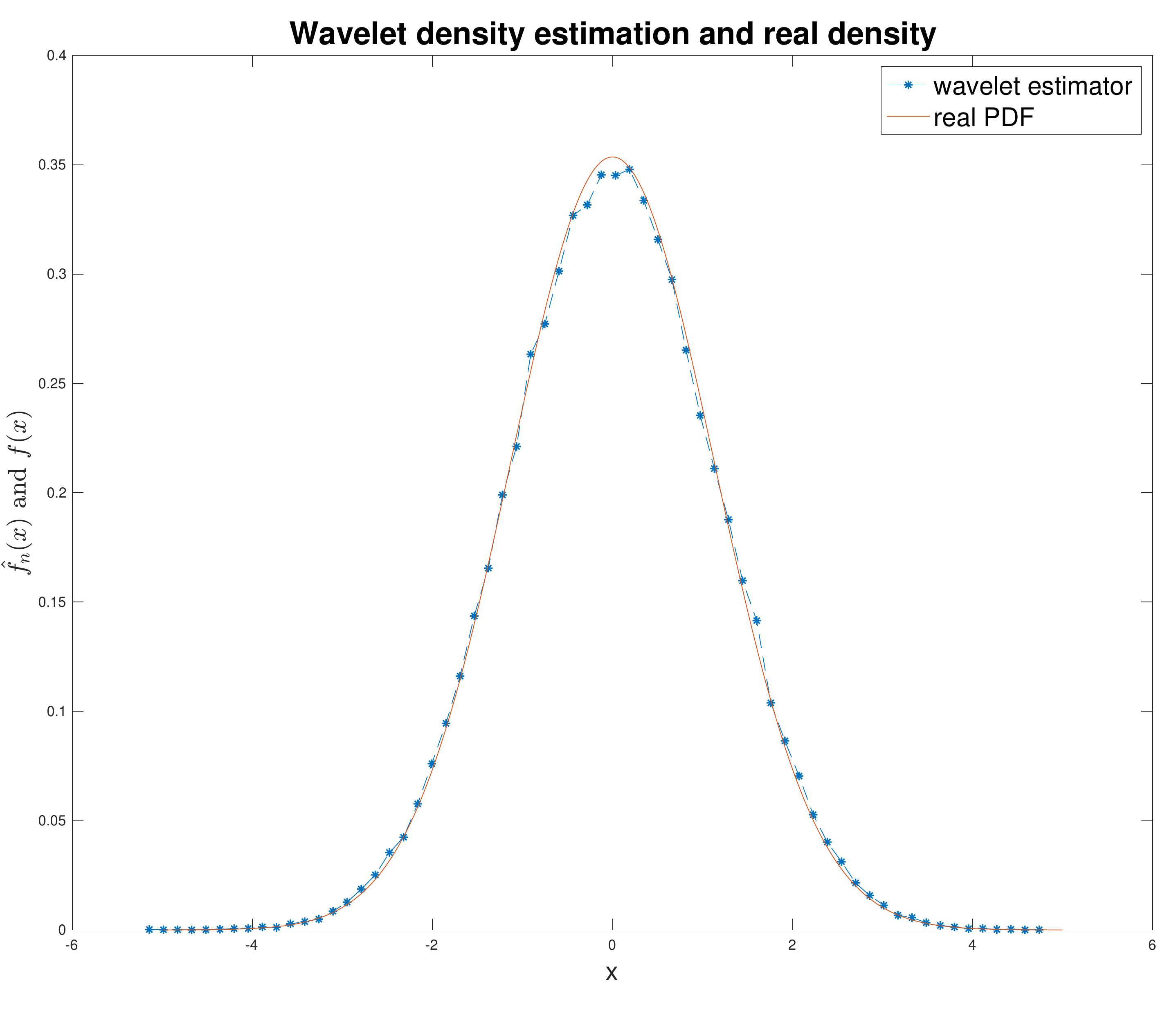}
\includegraphics[width=0.45\linewidth]{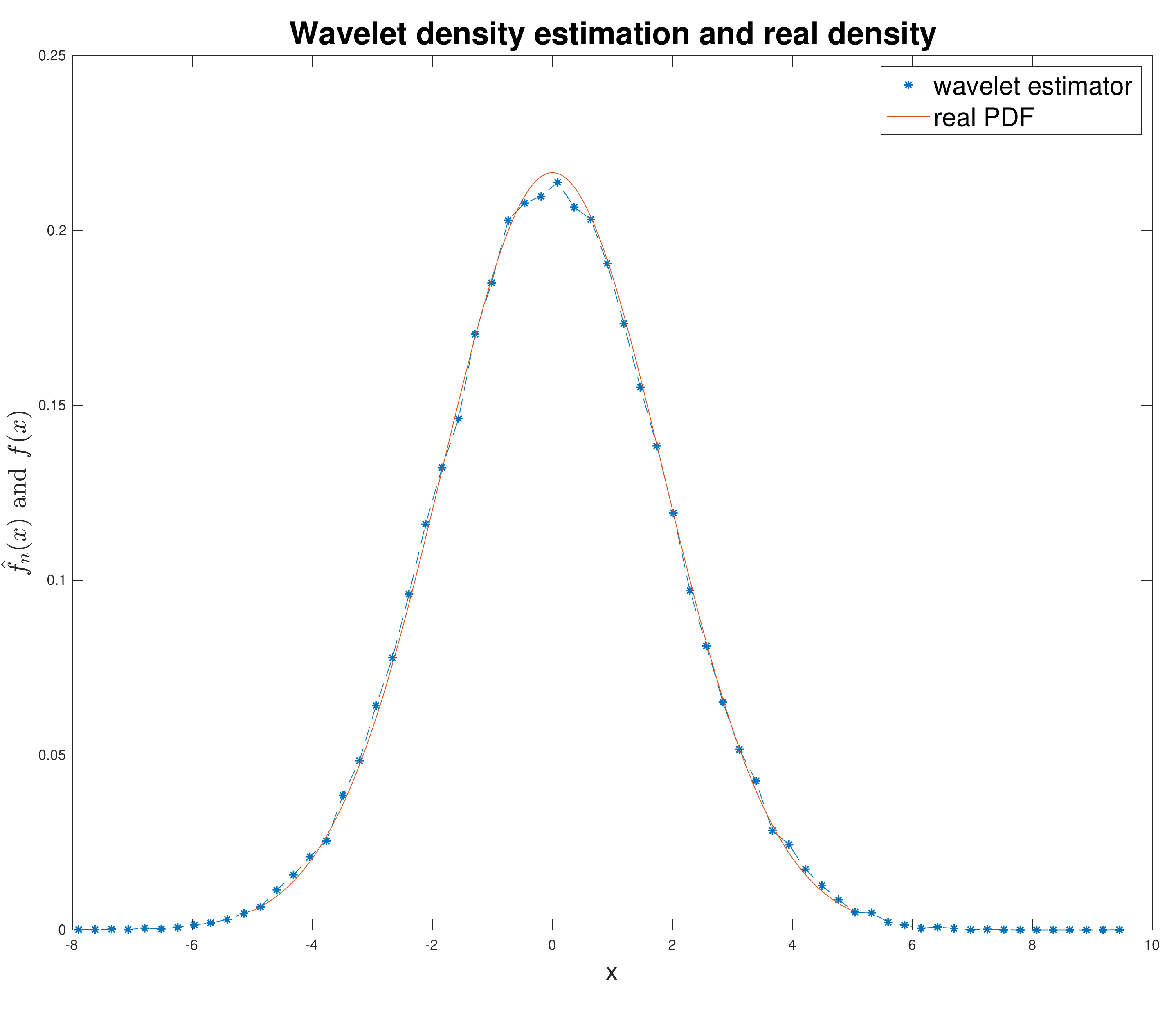}
\end{tabular}
\caption {Wavelet density estimation of simulated processes with Gaussian innovations and $d=-0.5$ (left), $-1.5$ (right). \label{fig:1}}
\end{figure}

\begin{figure}[H]
\centering
\begin{tabular}{cc}
\includegraphics[width=0.45\linewidth]{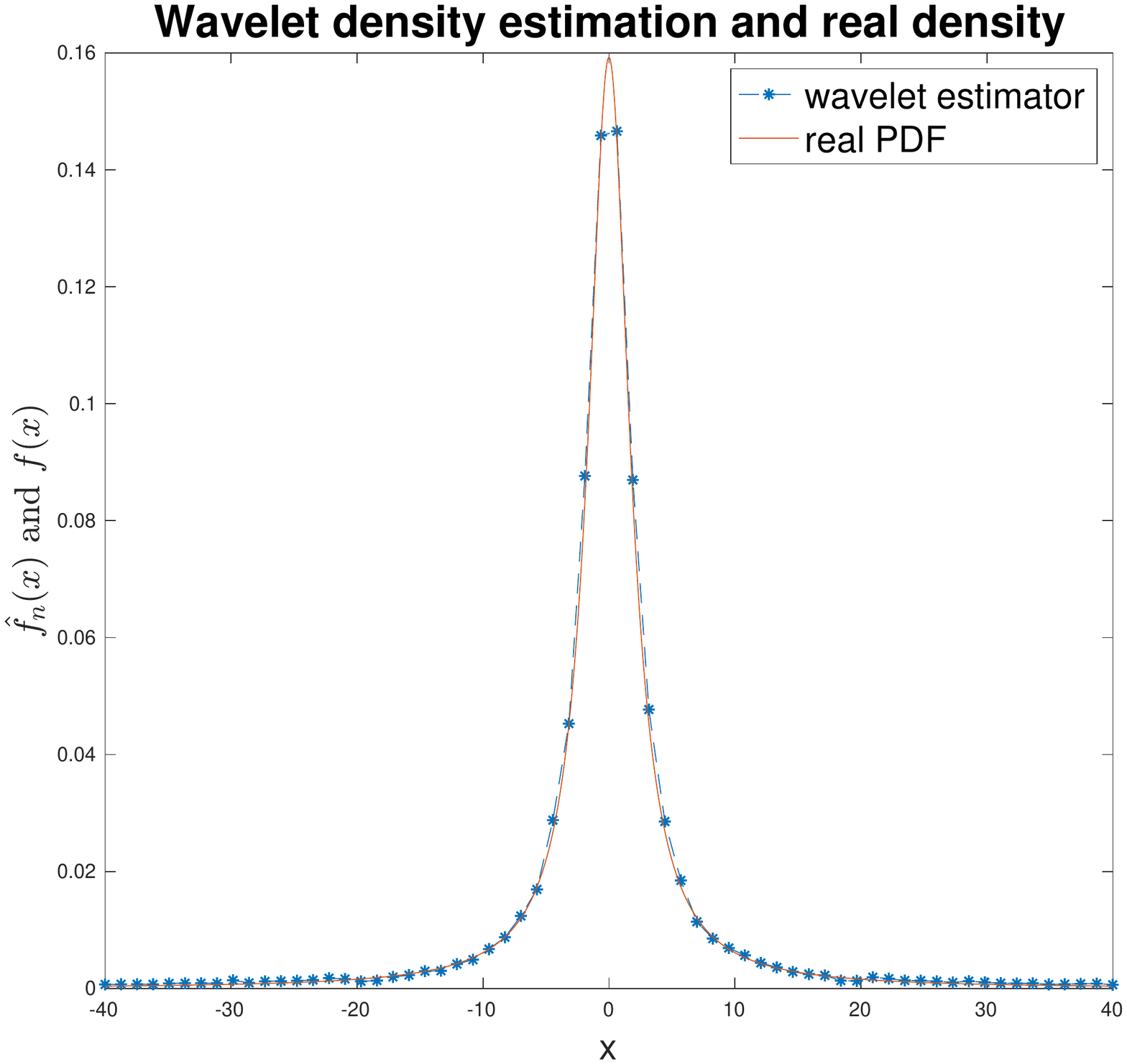}
\includegraphics[width=0.45\linewidth]{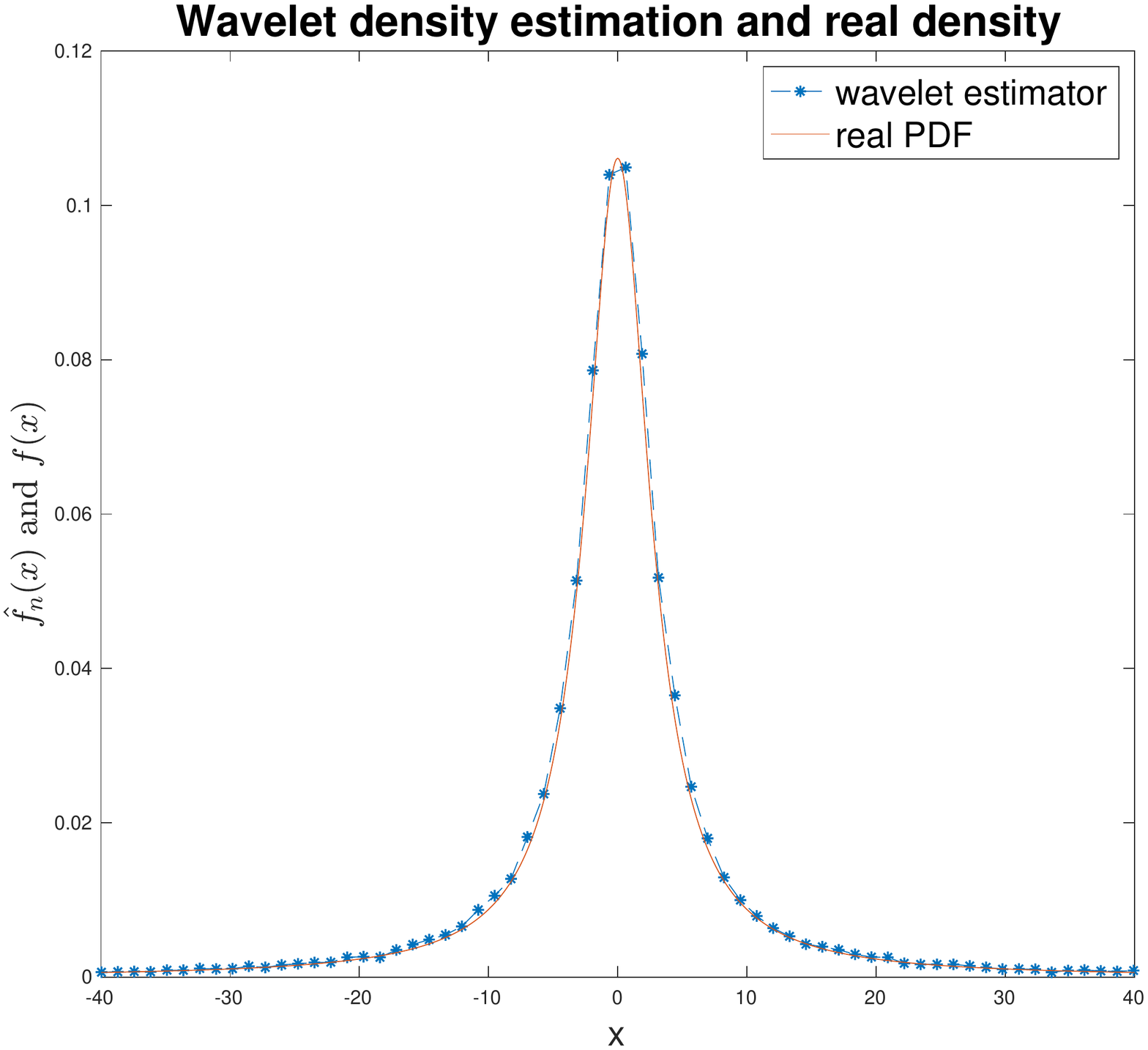}
\end{tabular}
\caption {Wavelet density estimation of simulated processes with Cauchy innovations and $d=-0.5$ (left), $-1.5$ (right). \label{fig:2}}
\end{figure}

\begin{figure}[H]
\centering
\begin{tabular}{cc}
\includegraphics[width=0.55\linewidth]{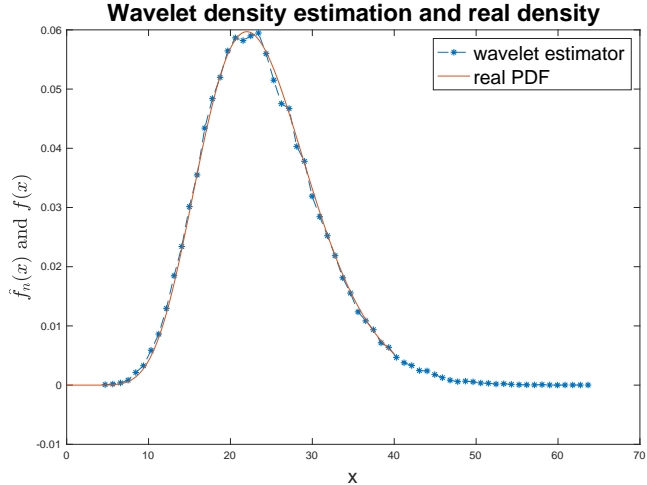}
\end{tabular}
\caption {Wavelet density estimation of simulated moving average process with $4$ non-zero coefficients. The innovations follow chi-squared distribution with $6$ degrees of freedom and the moving average process follows chi-squared distribution with $24$ degrees of freedom.\label{fig:3}}
\end{figure}

We also perform a simulation study for the moving average process $X_n=\sum\limits^{3}_{i=0}\varepsilon_{n-i}$ of order $4$ with coefficients $a_i=1$, $0\le i\le 3,$ and with innovations following chi-squared distribution with $6$ degrees of freedom. Then $X_n$ follows chi-squared distribution with $24$ degrees of freedom. Here we produce a MA$(4)$ process of length $n=2^{16}$ and, as it was done in the previous cases, to estimate the true density function, we apply the $1-$D wavelet estimator from MATLAB software with Daubechies wavelets of order $8$. Performance of the estimator in this case is presented in Figure \ref{fig:3} which again confirms the result of Theorem \ref{main}.  

\section{Proofs}\label{proof}
The orthonormality of functions $\varphi_{0k}(x)$ and $\psi_{jk}(x),$ $j\ge 0$, $k\in \Z$, together with the decompositions \eqref{1} and \eqref{estimator}, allows to represent the integrated mean squared error (IMSE) of the estimator $\hat{f}_n$ in the form
\begin{equation}\label{est}
\EX \int_\mathbb{R} \bigg[\hat{f}_n(x)-f(x)\bigg]^2dx=I_1+I_2+I_3,
\end{equation}
where
$$I_1=\sum\limits_{k\in \mathbb{Z}} \EX (\hat{\alpha}_{0k}-\alpha_{0k})^2, \quad \quad I_2=\sum\limits_{j=0}^{j_n}\sum\limits_{k\in \mathbb{Z}} \EX (\hat{\beta}_{jk}-\beta_{jk})^2$$
 and
$$I_3=\sum\limits_{j=j_n+1}^{\infty}\sum\limits_{k\in \mathbb{Z}} \beta^2_{jk}.$$

To prove Theorem \ref{main}, in Subsection \ref{estI123} we estimate each of $I_1$, $I_2$ and $I_3$ by representing them as sums of the values of Fourier transforms of integrable functions at integers and applying the formula of Poisson. In Subsection \ref{wavelet} we use the Assumption \ref{assu} to derive the properties of considered wavelets that will be used in the estimations of  $I_1$, $I_2$ and $I_3$, and the conditions of Poisson formula are verified in Subsection \ref{pf}.  We begin with the Subsection \ref{bounds} below which presents two auxiliary lemmas regarding the characteristic functions of innovations and linear processes that will be used in the later proofs. 

\subsection{Estimation of $\sum\limits_{1\leq i,  j \leq n}|H_{ij}(u,v)|$}\label{bounds}
Recall that for $i, j\in\N$ and for $\lambda, u, v\in\R$, we define $H(X_i)(\lambda)=e^{\iota \lambda X_i}-\phi(\lambda)$ and 
$H_{ij}(u,v)=\EX e^{\iota uX_i+\iota vX_j}-\phi(u)\phi(v)=\EX \big[H(X_i)(u)H(X_j)(v)\big]$, 
where $\phi(\lambda)=\E[e^{\iota \lambda  X_1}]$ is the characteristic function of the linear process $X_n=\sum\limits^{\infty}_{i=0} a_i\varepsilon_{n-i}$.    
To estimate the order of $\sum\limits_{1\leq i,  j \leq n}|H_{ij}(u,v)|$, we extend the projection method applied in Sang et al. (2018). In their Lemma 7.1, they applied the projection method to estimate the upper bound of the quantity $\sum_{1\leq i\neq j\leq n} \left| \E\left[H(X_i)(\lambda)\overline{H(X_j)(\lambda)}\right]\right|$, where $\overline{H(X_j)(\lambda)}$ is the conjugate of ${H(X_j)(\lambda)}$.  

For each $i\in\Z$, let  $\mathcal{F}_i$ be the $\sigma$-field generated by $\{\varepsilon_k:\, k\le i\}$. Given an integrable complex-valued random variable $Y$, we define the projection operator $\mathcal{P}_i$ as
\begin{align} \label{proj}
\mathcal{P}_i Y=\E[Y|\mathcal{F}_i]-\E[Y|\mathcal{F}_{i-1}].
\end{align}
It is easy to see that if $i\neq j$ then for any pair of integrable complex-valued random variables $Y$ and $W$ we have
\begin{align}
\E[\mathcal{P}_i Y\, \mathcal{P}_j W]=0.\label{projprod}
\end{align}

In the estimation of $\sum\limits_{1\leq i,  j \leq n}|H_{ij}(u,v)|$, we shall apply the condition (\ref{cha2}) of Theorem \ref{main}. As in the theorem this condition is imposed only for the case $\gamma\in (1/2, 1]$, in the following auxiliary lemma we show that the condition also holds for the case $\gamma\in (0,1/2]$.

\begin{lemma}\label{case}If the derivative of $\phi_{\varepsilon}$ has bounded modulus on $\mathbb{R}$, then \eqref{cha2} holds 
for all $\gamma\in (0,1/2],$ $\lambda\in \R$. 
\end{lemma}
\begin{proof}Note that for all $\lambda\in \mathbb{R}$,  $\E|e^{\iota \lambda \varepsilon_1}-\phi_{\varepsilon}(\lambda)|^{2}=1-|\phi_{\varepsilon}(\lambda)|^{2}\leq 1$.  In particular, if $|\lambda|\geq1$, then \eqref{cha2} holds for all $\gamma\in(0,1]$. For $|\lambda|<1$ and $\gamma\in (0,1/2]$ we have
\begin{align*}&\E|e^{\iota \lambda \varepsilon_1}-\phi_{\varepsilon}(\lambda)|^{2}=1-|\phi_{\varepsilon}(\lambda)|^{2}=\\
&(1+|\phi_{\varepsilon}(\lambda)|)(1-|\phi_{\varepsilon}(\lambda)|)\leq 2(|\phi_{\varepsilon}(0)|-|\phi_{\varepsilon}(\lambda)|)\leq2|\phi_{\varepsilon}(0)-\phi_{\varepsilon}(\lambda)|\leq c|\lambda|\leq c|\lambda|^{2\gamma},\end{align*} 
where the last inequality holds as $\gamma\in (0,1/2]$. 
\end{proof}

\begin{lemma} \label{lma1} 
Suppose $\sum\limits^{\infty}_{i=0}|a_i|^{\gamma}<\infty$ and \eqref{cha2} holds
for some $\gamma \in (0,1]$. 
Then there exists a positive constant $c$ such that
\begin{align} \label{lminq}
\sum\limits_{1\leq i \ne j \leq n}|H_{ij}(u,v)|\leq cn\big\{|uv|^{\gamma}|\phi_{\varepsilon}(u a_{0})\phi_{\varepsilon}(v a_{0})|+|u|^{\gamma}|\phi_{\varepsilon}(u a_{0})|+|v|^{\gamma}|\phi_{\varepsilon}(v a_{0})|\big\}
\end{align}
and
\begin{align}\label{equal}
&\left|H_{11}(u,v)\right|=\left|\E\left[H(X_1)(u)H(X_1)(v)\right]\right|\notag\\
&\leq c|uv|^{\gamma}|\phi_{\varepsilon}(u a_{0})\phi_{\varepsilon}(v a_{0})|+ c(|uv|^{\gamma}\wedge 1)\prod^{\infty}_{\ell=1} |\phi_{\varepsilon}(ua_{\ell}+va_{\ell})|.
\end{align}
In particular, if the condition \eqref{integral} is also satisfied, then 
\begin{align} \label{lminq1}
\sum\limits_{1\leq i, j \leq n}|H_{ij}(u,v)|\leq cn,
\end{align}
where $c>0$ is a constant that only depends on $a_0$ and $\gamma$.
\end{lemma}

\begin{proof}  
Using the definition of projection operator $\mathcal{P}_k$ in \eqref{proj}, applying \eqref{projprod}, the telescoping technique and the triangle inequality, for $i\le j$ we have 
\begin{align*}
\left|\E\left[H(X_i)(u)H(X_j)(v)\right]\right|
&\leq \sum^i_{k=-\infty}\left|\E\left[\mathcal{P}_{k}H(X_i)(u)\mathcal{P}_kH(X_j)(v)\right]\right| \\
&=\sum^i_{k=-\infty}\left|\E\left[\mathcal{P}_{0}H(X_{i-k})(u)\mathcal{P}_0H(X_{j-k})(v)\right]\right|\\
&=\sum^i_{k=-\infty}|\E\Big[e^{\iota\sum\limits^{\infty}_{\ell=1} (a_{\ell+i-k}u+a_{\ell+j-k}v) \varepsilon_{-\ell}} \Big] \Big[\prod^{i-k-1}_{\ell=0}\phi_{\varepsilon}(u a_{\ell})\prod^{j-k-1}_{\ell=0} \phi_{\varepsilon}(v a_{\ell}) \Big]\\
&\quad\times\E\left[(e^{\iota  u a_{i-k} \varepsilon_0}-\phi_{\varepsilon}(u a_{i-k}))(e^{\iota  v a_{j-k} \varepsilon_0}-\phi_{\varepsilon}(v a_{j-k}))\right]|.
\end{align*}

Here $\prod^{m}_{\ell=0}\phi_{\varepsilon}(u a_{\ell})=1$  if $m<0$. 
Next, we decompose the  sum $\sum^{i}_{k=-\infty}$ into $2$ sums corresponding to the cases $k<i$ and $k=i$. By Cauchy-Schwartz inequality, we have 
\begin{align*}
&\left|\E\left[H(X_i)(u)H(X_j)(v)\right]\right|\\
&\leq \sum^{i-1}_{k=-\infty} |\phi_{\varepsilon}(u a_{0})\phi_{\varepsilon}(v a_{0})| \sqrt{1-|\phi_{\varepsilon}(u a_{i-k})|^2}\sqrt{1-|\phi_{\varepsilon}(v a_{j-k})|^2}\\
&\quad+\prod^{j-i-1}_{\ell=0} |\phi_{\varepsilon}(v a_{\ell})| \prod^{\infty}_{\ell=1} |\phi_{\varepsilon}(ua_{\ell}+va_{\ell+j-i})| \sqrt{1-|\phi_{\varepsilon}(u a_0)|^2} \sqrt{1-|\phi_{\varepsilon}(v a_{j-i})|^2}.
\end{align*}
Therefore, using the conditions of the lemma, for the case $i<j$ we get
\begin{equation*}
\begin{split} 
\left|\E\left[H(X_i)(u)H(X_j)(v)\right]\right|
&\leq c|uv|^{\gamma}\sum^{i-1}_{k=-\infty}|\phi_{\varepsilon}(u a_{0})\phi_{\varepsilon}(v a_{0})| |a_{i-k}|^{\gamma}|a_{j-k}|^{\gamma}\\
&\qquad+ c|v|^{\gamma}|\phi_{\varepsilon}(v a_{0})| |a_{j-i}|^{\gamma}.
\end{split} 
\end{equation*} 

Hence,
\begin{equation*}
\begin{split} 
&\sum\limits_{1\leq i\leq n-1} \;\;\sum\limits_{i+1\leq j\leq n}|\EX \big[H(X_i)(u)H(X_j)(v)\big ]|\\
&\leq c\sum\limits_{1\leq i\leq n-1}\;\;\sum\limits_{i+1\leq j\leq n}|uv|^{\gamma}\sum^{i-1}_{k=-\infty}|\phi_{\varepsilon}(u a_{0})\phi_{\varepsilon}(v a_{0})| |a_{i-k}|^{\gamma}|a_{j-k}|^{\gamma}\\
&+c\sum\limits_{1\leq i\leq n-1}\;\;\sum\limits_{i+1\leq j\leq n}|v|^{\gamma}|\phi_{\varepsilon}(v a_{0})| |a_{j-i}|^{\gamma}\\
&\leq cn|uv|^{\gamma}|\phi_{\varepsilon}(u a_{0})\phi_{\varepsilon}(v a_{0})|+cn|v|^{\gamma}|\phi_{\varepsilon}(v a_{0})|.
\end{split} 
\end{equation*} 
Similarly, 
\begin{equation*}
\begin{split} 
&\sum\limits_{1\leq j\leq n-1} \;\;\sum\limits_{j+1\leq i\leq n}|\EX \big[H(X_i)(u)H(X_j)(v)\big ]|\\
&\leq cn|uv|^{\gamma}|\phi_{\varepsilon}(u a_{0})\phi_{\varepsilon}(v a_{0})|+cn|u|^{\gamma}|\phi_{\varepsilon}(u a_{0})|.
\end{split} 
\end{equation*} 
Thus,
\begin{equation*}
\begin{split} 
&\sum\limits_{1\leq i\neq j\leq n}|\EX \big[H(X_i)(u)H(X_j)(v)\big ]|\\
&=\bigg\{\sum\limits_{1\leq i\leq n-1} \;\sum\limits_{i+1\leq j\leq n}+\sum\limits_{1\leq j\leq n-1} \;\sum\limits_{j+1\leq i\leq n}\bigg\}|\EX \big[H(X_i)(u)H(X_j)(v)\big ]|\\
&\leq cn\big\{|uv|^{\gamma}|\phi_{\varepsilon}(u a_{0})\phi_{\varepsilon}(v a_{0})|+|u|^{\gamma}|\phi_{\varepsilon}(u a_{0})|+|v|^{\gamma}|\phi_{\varepsilon}(v a_{0})|\big\}
\end{split} 
\end{equation*} 
which proves \eqref{lminq}. To prove \eqref{equal} we note that when $i=j$ then 
\begin{align}
&\left|\E\left[H(X_i)(u)H(X_i)(v)\right]\right|\notag\\
&\leq \sum^{i-1}_{k=-\infty} |\phi_{\varepsilon}(u a_{0})||\phi_{\varepsilon}(v a_{0})| \sqrt{1-|\phi_{\varepsilon}(u a_{i-k})|^2}\sqrt{1-|\phi_{\varepsilon}(v a_{i-k})|^2}\label{t1}\\
&\quad+ \prod^{\infty}_{\ell=1} |\phi_{\varepsilon}(ua_{\ell}+va_{\ell})| \sqrt{1-|\phi_{\varepsilon}(u a_0)|^2} \sqrt{1-|\phi_{\varepsilon}(v a_{0})|^2}\label{t2}\\
&\leq c|uv|^{\gamma}|\phi_{\varepsilon}(u a_{0})\phi_{\varepsilon}(v a_{0})|+ c(|uv|^{\gamma}\wedge 1)\prod^{\infty}_{\ell=1} |\phi_{\varepsilon}(ua_{\ell}+va_{\ell})|.\notag
\end{align}
\end{proof}

\begin{remark}\label{others}
In the i.i.d. case, that is, in the case when $a_0$ is the only nonzero coefficient, we have that both $\sum\limits_{1\leq i\ne  j \leq n}|H_{ij}(u,v)|$ and \eqref{t1} vanish, while the term \eqref{t2} is equal to $c|uv|^{\gamma}$. Thus, in this case we get $\sum\limits_{1\leq i,  j \leq n}|H_{ij}(u,v)|=n|H_{11}(u,v)|\le cn(|uv|^{\gamma}\wedge 1)$. 
\end{remark}


\subsection{Properties of wavelets}\label{wavelet}
In this part we present the properties of wavelets that will be used in the following proofs. According to Assumption \ref{assu},  $\varphi$ has compact support and is twice continuously differentiable. Therefore (see, e.g., Stein and Shakarchi (2003),  pp. 132), its Fourier transform satisfies 
\begin{equation}\label{phi}
\begin{split} 
&|\hat{\varphi}(u)|\leq\frac{c}{1+u^2}.\\
\end{split} 
\end{equation}
Moreover, as $(\hat{\varphi})'=-\iota\widehat{x\varphi(x)}$ and $x\varphi(x)$ also has compact support and is twice continuously differentiable, then
\begin{equation}\label{phi1}
\begin{split} 
&|(\hat{\varphi}(u))'|\leq\frac{c}{1+u^2}.\\
\end{split} 
\end{equation}
Similarly, for the Fourier transform of $\psi$ we have 
\begin{equation}\label{sumpsi}
\begin{split} 
&|\hat{\psi}(u)|\leq\frac{c}{1+u^2}\\
\end{split} 
\end{equation}
and
\begin{equation}\label{psi1}
\begin{split} 
&|(\hat{\psi}(u))'|\leq\frac{c}{1+u^2}.\\
\end{split} 
\end{equation}
Also, as the mother wavelet function of the considered Daubechies wavelets has $\lceil{M\beta}\rceil$ vanishing moments $\int_\mathbb{R} x^r\psi(x)dx=0,$ then $(\hat{\psi})^{(r)}(0)=0,$ $r=0,...,\lceil{M\beta}\rceil-1$. Since the derivative of $\hat{\psi}$ of order$\lceil{M\beta}\rceil$ is bounded: $|(\hat{\psi})^{(\lceil{M\beta}\rceil)}(u)|\le \int_\R t^{\lceil{M\beta}\rceil}|\psi(t)|dt<\infty$, then the function $\hat{\psi}(s)/s^{\lceil{M\beta}\rceil}$ is bounded on $\mathbb{R}$. Here we shall use the L'Hopital's Rule to obtain the boundedness of $\hat{\psi}(s)/s^{\lceil{M\beta}\rceil}$ around $0$. In particular, the function $\hat{\psi}(s)/s^{M\beta}$ is also bounded on $\mathbb{R}$.
\subsection{Poisson summation formula}\label{pf}
To prove Theorem \ref{main} we will consider expressions of the form 
$$\sum\limits_{k\in\mathbb{Z}}\int_\mathbb{R}\int_\mathbb{R} e^{\iota k(u+v)}g(u,v)dudv,$$
where $g:\mathbb{R}\times\mathbb{R}\to\mathbb{R}$ is some integrable function.
Defining  
\begin{equation}\label{h}h(t) = \int_\mathbb{R} g(t+s,-s) ds,\end{equation}
we can rewrite the above sum as
$$\sum_{k\in\mathbb{Z}} \int_\mathbb{R} e^{\iota  k t} h(t) dt=\sum_{k\in\mathbb{Z}}\hat{h}(k),$$
and the latter sum of the values of Fourier transform of the function $h(t)$ at integers can be calculated using the following Poisson summation formula (see Zygmund (1966), Ch. II, \S 13):
\begin{theorem}[Poisson summation formula]\label{formula} If $h\in L^1(\mathbb{R})\cap BV(\mathbb{R})$ is an absolutely integrable function of bounded variation and $2h(x)=h(x+0)+h(x-0)$ for all $x\in\mathbb{R}$, then 
$$\sum_{k\in\mathbb{Z}}\hat{h}(k)=2\pi \sum\limits_{k\in\mathbb{Z}}h(2k\pi).$$
\end{theorem}
In particular, if the function $h$ defined by \eqref{h} satisfies the conditions of the above theorem, then
\begin{align}\label{Poisson}
\sum\limits_{k\in\mathbb{Z}}\int_\mathbb{R}\int_\mathbb{R} e^{\iota k(u+v)}g(u,v)dudv=2\pi \sum\limits_{k\in\mathbb{Z}}h(2k\pi).
\end{align}
Throughout the proof of Theorem \ref{main} we are going to apply the Poisson summation formula to functions 
$$h_1(t)=\sum\limits_{1\leq i, m \leq n}\int_\mathbb{R} \hat{\varphi}(t+s)\hat{\varphi}(-s)H_{im}(t+s,-s) ds,$$
$$h_2(t)=\sum\limits_{1\leq i, m \leq n}\int_\mathbb{R} \hat{\psi}(t+s)\hat{\psi}(-s)H_{im}(2^j(t+s),-2^js) ds$$
and
$$h_3(t)=\int_\mathbb{R} \hat{\psi}(t+s)\hat{\psi}(-s) \phi(2^{j}(t+s))\phi(-2^{j}s)ds,$$
where $\varphi$ and $\psi$ are, respectively, the scaling function and the mother wavelet discussed in the previous subsection, $\phi$ is the characteristic function of linear process $X_n$ and $j\in\mathbb{N}$ is some number.
\begin{lemma}\label{Poisson1}Under the conditions of Theorem \ref{main}, the functions $h_1, h_2$ and $h_3$ satisfy the conditions of Theorem \ref{formula}.\end{lemma}
\begin{proof}
As $H_{ij}(u,v):=\EX e^{\iota uX_i+\iota vX_j}-\phi(u)\phi(v),$ then $|H_{ij}(u,v)|\leq2$ for all $i, j\in\mathbb{N}$, which together with \eqref{phi} implies that $h_1\in L^1(\mathbb{R})$. Let us now show that $\int_\mathbb{R}|h_1'(t)|dt<\infty$, which will imply (see Baernstein (2019), page 125) that $h_1\in BV(\mathbb{R})$. To bound the function $h_1'(t)$ we first estimate $\bigg|\frac{\partial}{\partial u} H_{ij}(u,v)\bigg|$. As $X_i=\sum\limits^{\infty}_{p=0} a_p\varepsilon_{i-p}=\sum\limits^{\infty}_{m=-i} a_{i+m}\varepsilon_{-m}$ and $X_j=\sum\limits^{\infty}_{m=-j} a_{j+m}\varepsilon_{-m}$, then for given $u,v\in\mathbb{R}$ and for $i\leq j$ we have that
$$uX_i+vX_j=\sum\limits^{\infty}_{m=-i} ua_{i+m}\varepsilon_{-m}+\sum\limits^{\infty}_{m=-j} va_{j+m}\varepsilon_{-m}
=\sum\limits^{-i-1}_{m=-j} va_{j+m}\varepsilon_{-m}+\sum\limits^{\infty}_{m=-i} (ua_{i+m}+va_{j+m})\varepsilon_{-m}$$
and, therefore, 
$$
E_{ij}(u,v):=\EX e^{iuX_i+ivX_j}=\prod\limits^{-i-1}_{m=-j} \phi_{\varepsilon}(va_{j+m})\prod\limits^{\infty}_{m=-i} \phi_{\varepsilon} (ua_{i+m}+va_{j+m})$$ (as usual, sums and products with negative number of summands and factors are assumed to be, respectively, $0$ and $1$). Thus,
\begin{equation}\label{derEij}
\begin{split} 
&\frac{\partial}{\partial u}E_{ij}(u,v)=\sum\limits^{\infty}_{m=-i}a_{i+m} \phi'_{\varepsilon}(ua_{i+m}+va_{j+m})\frac{E_{ij}(u,v)}{\phi_{\varepsilon}(ua_{i+m}+va_{j+m})}.
\end{split} 
\end{equation}
 Also, as  
$$\phi(u)=\prod\limits^{\infty}_{i=0}\phi_{\varepsilon}(ua_i),$$
then (see Giraitis, Koul and Surgailis (1996),  page 322)
$$
\phi'(u)=\sum\limits_{i=0}^\infty a_i\phi'_{\varepsilon}(ua_i)\prod\limits_{\substack{j\geq 0 \\ j\neq i}}\phi_{\varepsilon}(ua_j), \quad u\in\mathbb{R},
$$ 
and, therefore, 
$$\frac{\partial}{\partial u}\phi(u)\phi(v)=\sum\limits^{\infty}_{m=-i}a_{i+m} \phi'_{\varepsilon}(ua_{i+m})\frac{\phi(u)\phi(v)}{\phi_{\varepsilon}(ua_{i+m})}.$$ Hence, using the boundedness of $\phi'_{\varepsilon}$, we get
\begin{align}
&\sup\limits_{u,v\in\mathbb{R}}\bigg|\frac{\partial}{\partial u} H_{ij}(u,v)\bigg|=\sup\limits_{u,v\in\mathbb{R}}\bigg|\frac{\partial}{\partial u}E_{ij}(u,v)-\frac{\partial}{\partial u}\phi(u)\phi(v)\bigg|\notag\\
&=\sup\limits_{u,v\in\mathbb{R}}\bigg|\sum\limits^{\infty}_{m=-i}a_{i+m} \phi'_{\varepsilon}(ua_{i+m}+va_{j+m})\frac{E_{ij}(u,v)}{\phi_{\varepsilon}(ua_{i+m}+va_{j+m})}-\sum\limits^{\infty}_{m=-i}a_{i+m} \phi'_{\varepsilon}(ua_{i+m})\frac{\phi(u)\phi(v)}{\phi_{\varepsilon}(ua_{i+m})}\bigg|\notag\\
&\leq c\sum\limits^{\infty}_{m=-i}|a_{i+m}|<\infty\notag.
\end{align}
Therefore, using also \eqref{phi1}, we have
$$\bigg|\frac{\partial}{\partial t}\bigg[\hat{\varphi}(t+s)\hat{\varphi}(-s)H_{ij}(t+s,-s)\bigg]\bigg|\leq c|\hat{\varphi}(-s)|,$$
and since $\int_\mathbb{R}|\hat{\varphi}(-s)|ds<\infty$, we get  by dominated convergence theorem, that 

\begin{align*}&\int_\mathbb{R}|h_1'(t)|dt \\
&\leq \sum\limits_{1\leq i, j \leq n}\int_\mathbb{R}\bigg|\int_\mathbb{R}\frac{\partial}{\partial t}\bigg[\hat{\varphi}(t+s)\hat{\varphi}(-s)H_{ij}(t+s,-s)\bigg]ds\bigg|dt\leq cn^2\int_\mathbb{R}\int_\mathbb{R}\frac{1}{1+(t+s)^2}\frac{1}{1+s^2}dsdt<\infty,\end{align*}
where we also used the inequalities \eqref{phi} and \eqref{phi1}.  Proofs for the functions $h_2$ and $h_3$ are identical.
\end{proof}

\subsection{Estimation of $I_1$, $I_2$ and $I_3$}\label{estI123}
\begin{lemma} \label{lma2} 
Under the conditions of Theorem \ref{main}

\begin{align*} \label{I1res}
|I_1|\leq \frac{c}{n}.
\end{align*}
\end{lemma}
\begin{proof}
We have that
\begin{equation*}
\begin{split} 
&I_1=\sum\limits_{k\in\mathbb{Z}} \EX (\hat{\alpha}_{0k}-\alpha_{0k})^2=\sum\limits_{k\in\mathbb{Z}} \EX\bigg[\frac{1}{n}\sum\limits_{i=1}^n\varphi_{0k}(X_i)-\EX \varphi_{0k}(X_1)\bigg]^2\\
&=\frac{1}{4\pi^2}\sum\limits_{k\in\mathbb{Z}} \EX\bigg[\frac{1}{n}\sum\limits_{i=1}^n\int_\mathbb{R}\hat{\varphi}_{0k}(u)e^{\iota uX_i}du-\EX \int_\mathbb{R}\hat{\varphi}_{0k}(u)e^{\iota uX_1}du\bigg]^2\\
&=\frac{1}{4\pi^2}\sum\limits_{k\in\mathbb{Z}} \EX\bigg[\frac{1}{n}\sum\limits_{i=1}^n\int_\mathbb{R}\hat{\varphi}_{0k}(u)e^{\iota uX_i}du- \frac{1}{n}\sum\limits_{i=1}^n\int_\mathbb{R}\hat{\varphi}_{0k}(u)\phi(u)du\bigg]^2\\
&=\frac{1}{4\pi^2n^2}\sum\limits_{k\in\mathbb{Z}} \EX\bigg[\sum\limits_{i=1}^n\int_\mathbb{R}\hat{\varphi}_{0k}(u)[e^{\iota uX_i}-\phi(u)]du\bigg]^2\\
&=\frac{1}{4\pi^2n^2}\sum\limits_{k\in\mathbb{Z}} \EX\bigg[\int_\mathbb{R}\hat{\varphi}(u)e^{\iota uk}\bigg(\sum\limits_{i=1}^n[e^{\iota uX_i}-\phi(u)]\bigg)du\bigg]^2\\
&=\frac{1}{4\pi^2n^2}\sum\limits_{k\in\mathbb{Z}} \EX\Bigg[\int_\mathbb{R}\int_\mathbb{R}\hat{\varphi}(u)\hat{\varphi}(v)e^{\iota (u+v)k}\bigg(\sum\limits_{i=1}^n[e^{\iota uX_i}-\phi(u)]\bigg)\bigg(\sum\limits_{i=1}^n[e^{\iota vX_i}-\phi(v)]\bigg)dudv\Bigg]\\
&=\frac{1}{4\pi^2n^2}\sum\limits_{k\in\mathbb{Z}} \Bigg[\int_\mathbb{R}\int_\mathbb{R}\hat{\varphi}(u)\hat{\varphi}(v)e^{\iota (u+v)k}\sum\limits_{1\leq i, j \leq n}H_{ij}(u,v)dudv\Bigg].\\
\end{split} 
\end{equation*}
Let 
\begin{align*}
g_1(u,v)=\hat{\varphi}(u)\hat{\varphi}(v)\sum\limits_{1\leq i, j \leq n}H_{ij}(u,v). 
\end{align*}
Then 
\begin{align*}
h_1(t)=\int_\mathbb{R} g_1(t+s,-s)ds=\int_\mathbb{R} \hat{\varphi}(t+s)\hat{\varphi}(-s)\sum\limits_{1\leq i, j \leq n}H_{ij}(t+s,-s) ds
\end{align*}
and by Theorem \ref{formula} and Lemma \ref{Poisson1} we have 
\begin{align*}
I_1&=\frac{1}{2\pi n^2}\sum\limits_{k\in\mathbb{Z}}h_1(2k\pi)\\
&=\frac{1}{2\pi n^2}\sum\limits_{k\in\mathbb{Z}}\int_\mathbb{R} \hat{\varphi}(2k\pi+s)\hat{\varphi}(-s)\sum\limits_{1\leq i, j \leq n}H_{ij}(2k\pi+s,-s) ds.
\end{align*}
Hence, applying \eqref{lminq1} and \eqref{phi} we get
\begin{align*}
|I_1|&\le \frac{1}{2\pi n^2}\sum\limits_{k\in\mathbb{Z}}\int_\mathbb{R} |\hat{\varphi}(2k\pi+s)||\hat{\varphi}(-s)|\sum\limits_{1\leq i, j \leq n}|H_{ij}(2k\pi+s,-s)| ds\\
&\le  \frac{c}{n}\sum\limits_{k\in\mathbb{Z}}\int_\mathbb{R} |\hat{\varphi}(2k\pi+s)||\hat{\varphi}(s)| ds\\
&\le \frac{c}{n}\int_\mathbb{R} \sum\limits_{k\in\mathbb{Z}}\frac{1}{1+(2k\pi+s)^2}\frac{1}{1+s^2} ds\\
&\le \frac{c}{n}\int_\mathbb{R}\frac{1}{1+s^2} ds\le \frac{c}{n}.
\end{align*}
\end{proof}

\begin{lemma} \label{lma3} 
Under the conditions of Theorem \ref{main},
for any $j_n>0$,
\begin{align*} 
|I_2|\leq \frac{c2^{j_n}}{n}.
\end{align*}
\end{lemma}
\begin{proof}
We have that
\begin{align*}
\begin{split} 
&I_2=\sum_{j=0}^{j_n}\sum\limits_{k\in\mathbb{Z}} \EX (\hat{\beta}_{jk}-\beta_{jk})^2=\sum_{j=0}^{j_n}\sum\limits_{k\in\mathbb{Z}} \EX\bigg[\frac{1}{n}\sum\limits_{i=1}^n\psi_{jk}(X_i)-\EX \psi_{jk}(X_1)\bigg]^2\\
&=\frac{1}{4\pi^2}\sum_{j=0}^{j_n}\sum\limits_{k\in\mathbb{Z}} \EX\bigg[\frac{1}{n}\sum\limits_{i=1}^n\int_\mathbb{R}\hat{\psi}_{jk}(u)e^{\iota uX_i}du-\EX \int_\mathbb{R}\hat{\psi}_{jk}(u)e^{\iota uX_1}du\bigg]^2\\
&=\frac{1}{4\pi^2}\sum_{j=0}^{j_n}\sum\limits_{k\in\mathbb{Z}} \EX\bigg[\frac{1}{n}\sum\limits_{i=1}^n\int_\mathbb{R}\hat{\psi}_{jk}(u)e^{\iota uX_i}du- \frac{1}{n}\sum\limits_{i=1}^n\int_\mathbb{R}\hat{\psi}_{jk}(u)\phi(u)du\bigg]^2\\
&=\frac{1}{4\pi^2n^2}\sum_{j=0}^{j_n}\sum\limits_{k\in\mathbb{Z}} \EX\bigg[\int_\mathbb{R}\hat{\psi}_{jk}(u)\sum\limits_{i=1}^n[e^{\iota uX_i}-\phi(u)]du\bigg]^2\\
&=\frac{1}{4\pi^2n^2}\sum_{j=0}^{j_n}\sum\limits_{k\in\mathbb{Z}} \EX\int_\mathbb{R}\int_\mathbb{R}\hat{\psi}_{jk}(u)\hat{\psi}_{jk}(v)\sum_{1\le i, m\le n}[e^{\iota uX_i}-\phi(u)][e^{\iota vX_m}-\phi(v)]dudv\\
&=\frac{2^{-j}}{4\pi^2n^2}\sum_{j=0}^{j_n}\sum\limits_{k\in\mathbb{Z}}\int_\mathbb{R}\int_\mathbb{R} e^{\iota ku2^{-j}}\hat{\psi}(u2^{-j})e^{\iota kv2^{-j}}\hat{\psi}(v2^{-j})\sum\limits_{1\leq i, m \leq n}H_{im}(u,v)dudv\\
&=\frac{2^{j}}{4\pi^2n^2}\sum_{j=0}^{j_n}\sum\limits_{k\in\mathbb{Z}} \Bigg[\int_\mathbb{R}\int_\mathbb{R}\hat{\psi}(u)\hat{\psi}(v)e^{\iota (u+v)k}\sum\limits_{1\leq i, m \leq n}H_{im}(2^{j}u,2^{j}v)dudv\Bigg]:=\sum_{j=0}^{j_n}I_{2,j}.\\
\end{split} 
\end{align*}

Let 
\begin{align*}
g_2(u,v)=\hat{\psi}(u)\hat{\psi}(v)\sum\limits_{1\leq i, m \leq n}H_{im}(2^ju,2^jv). 
\end{align*}
Then 
\begin{align*}
h_2(t)=\int_\mathbb{R} g_2(t+s,-s)ds=\int_\mathbb{R} \hat{\psi}(t+s)\hat{\psi}(-s)\sum\limits_{1\leq i, m \leq n}H_{im}(2^j(t+s),-2^js) ds.
\end{align*}
By the Poisson formula \eqref{Poisson} we have
\begin{align*}
I_{2,j}=\frac{2^j}{2\pi n^2}\sum\limits_{k\in\mathbb{Z}}h_2(2k\pi)=\frac{2^j}{2\pi n^2}\sum\limits_{k\in\mathbb{Z}}\int_\mathbb{R} \hat{\psi}(2k\pi+s)\hat{\psi}(-s)\sum\limits_{1\leq i, m \leq n}H_{im}(2^j(2k\pi+s),-2^js) ds.
\end{align*}
Hence, applying \eqref{lminq1} and \eqref{sumpsi}, we get 
\begin{align*}
|I_{2,j}|&\le \frac{2^j}{2\pi n^2}\sum\limits_{k\in\mathbb{Z}}\int_\mathbb{R} |\hat{\psi}(2k\pi+s)||\hat{\psi}(-s)|\sum\limits_{1\leq i, m \leq n}|H_{im}(2^j(2k\pi+s),-2^js)| ds\\
&\le  \frac{c2^j}{n}\sum\limits_{k\in\mathbb{Z}}\int_\mathbb{R} |\hat{\psi}(2k\pi+s)||\hat{\psi}(s)|ds\le \frac{c2^j}{n}\int_\mathbb{R} \sum\limits_{k\in\mathbb{Z}}\frac{1}{1+(2k\pi+s)^2}|\hat{\psi}(s)|ds\\
&\le \frac{c2^j}{n}\int_\mathbb{R}|\hat{\psi}(s)| ds\le \frac{c2^j}{n}.
\end{align*}

Thus,
$$|I_2|\le |\sum_{j=0}^{j_n}|I_{2,j}|\leq  \frac{c2^{j_n}}{n}.$$\end{proof}

\begin{lemma} \label{lma4} 
Under the conditions of Theorem \ref{main}, for any $j_n>0$,
\begin{align*} \label{I3res}
|I_3|\leq c2^{-2j_nM\beta}.
\end{align*}
\end{lemma}
\begin{proof}

Recall that  $\phi$ is the characteristic function of $X_1$.
We then  have that
\begin{equation*}
\begin{split} 
&I_{3,j}:=\sum\limits_{k\in\mathbb{Z}}\beta^2_{jk}=\sum\limits_{k\in\mathbb{Z}} \bigg[\EX \psi_{jk}(X_1)\bigg]^2\\
&=\frac{1}{4\pi^2}\sum\limits_{k\in\mathbb{Z}} \bigg[\EX \int_\mathbb{R}\hat{\psi}_{jk}(u)e^{\iota uX_1}du\bigg]^2\\
&=\frac{1}{4\pi^2}\sum\limits_{k\in\mathbb{Z}} \bigg[\int_\mathbb{R}\hat{\psi}_{jk}(u)\phi(u)du\bigg]^2\\
&=\frac{2^{j}}{4\pi^2}\sum\limits_{k\in\mathbb{Z}} \bigg[\int_\mathbb{R} e^{\iota ku}\hat{\psi}(u)\phi(2^{j}u)du\bigg]^2\\
&=\frac{2^{j}}{4\pi^2}\sum\limits_{k\in\mathbb{Z}} \Bigg[\int_\mathbb{R}\int_\mathbb{R} e^{\iota (u+v)k}\hat{\psi}(u)\hat{\psi}(v)\phi(2^{j}u)\phi(2^{j}v)dudv\Bigg].\\
\end{split} 
\end{equation*}
Let 
\begin{align*}
g_3(u,v)=\hat{\psi}(u)\hat{\psi}(v)\phi(2^{j}u)\phi(2^{j}v). 
\end{align*}
Then 
\begin{align*}
h_3(t)=\int_\mathbb{R} g_3(t+s,-s)ds=\int_\mathbb{R} \hat{\psi}(t+s)\hat{\psi}(-s) \phi(2^{j}(t+s))\phi(-2^{j}s)ds.
\end{align*}
Hence, by the Poisson formula \eqref{Poisson},
\begin{align*}
I_{3,j}&=\frac{2^j}{2\pi}\sum\limits_{k\in\mathbb{Z}}h_3(2k\pi)=\frac{2^j}{2\pi}\sum\limits_{k\in\mathbb{Z}}\int_\mathbb{R} \hat{\psi}(2k\pi+s)\hat{\psi}(-s) \phi(2^{j}(2k\pi+s))\phi(-2^{j}s)ds\\
& =\frac{2^j}{2\pi}\int_\mathbb{R} \hat{\psi}(s)\hat{\psi}(-s) \phi(2^{j}s)\phi(-2^{j}s)ds+\frac{2^j}{2\pi}\sum\limits_{k\in\mathbb{Z}\setminus\{0\}}\int_{-1}^1 \hat{\psi}(2k\pi+s)\hat{\psi}(-s) \phi(2^{j}(2k\pi+s))\phi(-2^{j}s)ds\\
&+\frac{2^j}{2\pi}\sum\limits_{k\in\mathbb{Z}\setminus\{0\}}\int_{|s|>1} \hat{\psi}(2k\pi+s)\hat{\psi}(-s) \phi(2^{j}(2k\pi+s))\phi(-2^{j}s)ds:=I'_{3,j}+I''_{3,j}+I'''_{3,j}.
\end{align*} 
As $\phi(u)=\prod\limits^{\infty}_{i=0}\phi_{\varepsilon}( a_iu)$ and there are at least $M$ nonzero coefficients among $a_i$,  $i\in\mathbb{N}_0$, then it follows from \eqref{integral} that 
$$\int_\mathbb{R}|u^{M\beta}\phi(u)|^2du<\infty.$$ 
Using also the boundedness of $\hat{\psi}(s)/s^{M\beta}$, we get 
\begin{align*}
|I'_{3,j}|\leq c2^j\int_\mathbb{R} |\hat{\psi}(s)\hat{\psi}(-s) \phi(2^{j}s)\phi(-2^{j}s)|ds\leq  c2^j\int_\mathbb{R} |s^{2M\beta}\phi(2^{j}s)\phi(-2^{j}s)|ds\leq 2^{-2jM\beta}.
\end{align*}
For $k\in\mathbb{Z}\setminus\{0\}$ and $s\in[-1,1]$ we have $|2k\pi+s|>1$, and, therefore, the condition \eqref{integral} implies that  $|\phi(2^{j}(2k\pi+s))|\leq c2^{-jM\beta}$. Hence,  
\begin{align*}
&|I''_{3,j}|\leq c2^j\sum\limits_{k\in\mathbb{Z}\setminus\{0\}}\int_{-1}^1 |\hat{\psi}(2k\pi+s)\hat{\psi}(-s) \phi(2^{j}(2k\pi+s))\phi(-2^{j}s)|ds\\
&\leq c2^{j-jM\beta}\sum\limits_{k\in\mathbb{Z}\setminus\{0\}}\int_{-1}^1 |\hat{\psi}(2k\pi+s)\hat{\psi}(-s)\phi(-2^{j}s)|ds\\
&\leq c2^{j-jM\beta}\int_{-1}^1\sum\limits_{k\in\mathbb{Z}\setminus\{0\}}\frac{1}{1+(2k\pi+s)^2}|\hat{\psi}(-s)\phi(-2^{j}s)|ds\\
&\leq c2^{j-jM\beta}\int_{-1}^1|\hat{\psi}(-s)\phi(-2^{j}s)|ds\leq c2^{j-jM\beta}\int_{-1}^1|s^{M\beta}\phi(-2^{j}s)|ds\leq  c2^{-2jM\beta}.
\end{align*}
Similarly, for $|s|>1$ we have $|\phi(-2^{j}s)|\leq c2^{-jM\beta}$, and, therefore,
\begin{align*}
&|I'''_{3,j}|\leq c2^j\sum\limits_{k\in\mathbb{Z}\setminus\{0\}}\int_{|s|>1} |\hat{\psi}(2k\pi+s)\hat{\psi}(-s) \phi(2^{j}(2k\pi+s))\phi(-2^{j}s)|ds\\
&\leq c2^{j-jM\beta}\sum\limits_{k\in\mathbb{Z}\setminus\{0\}}\int_{|s|>1} |\hat{\psi}(2k\pi+s)\hat{\psi}(-s)\phi(2^{j}(2k\pi+s))|ds\\
&\leq c2^{j-jM\beta}\sum\limits_{k\in\mathbb{Z}\setminus\{0\}}\int_{\mathbb{R}} |\hat{\psi}(2k\pi+s)\hat{\psi}(-s)\phi(2^{j}(2k\pi+s))|ds\\
&= c2^{j-jM\beta}\sum\limits_{k\in\mathbb{Z}\setminus\{0\}}\int_{\mathbb{R}} |\hat{\psi}(u)\hat{\psi}(-u+2k\pi)\phi(2^{j}u)|du\\
&\leq c2^{j-jM\beta}\int_\mathbb{R}|\hat{\psi}(u)\phi(2^{j}u)|du\leq c2^{j-jM\beta}\int_\mathbb{R}|u^{M\beta}\phi(2^{j}u)|du\leq  2^{-2jM\beta}\leq c2^{-2jM\beta}.
\end{align*}
Hence,
$$|I_3|=\bigg|\sum\limits_{j=j_n+1}^{\infty}I_{3,j}\bigg|\leq c2^{-2j_nM\beta}.$$\end{proof}

\noindent\textbf{Proof of Theorem \ref{main}}.
Putting together the results of Lemmas \ref{lma2}, \ref{lma3} and \ref{lma4} and taking $j_n=\lceil\frac{\log_2n}{2M\beta+1}\rceil$ proves the Theorem \ref{main}.

\bigskip

\textbf{Acknowledgement} We thank Johannes Schmidt-Hieber and Fangjun Xu for very helpful comments. We are also grateful to the reviewers for many constructive remarks and comments that helped to improve the manuscript.  The research of Aleksandr Beknazaryan and Peter Adamic is partially supported by the Natural Sciences and Engineering Research Council of Canada under Grant RGPIN-2017-05595. The research of Hailin Sang is partially supported by the Simons Foundation Grant 586789, USA.

\end{document}